\documentstyle{articleplus}

\title{Extremal Configurations of Hinge Structures}
\author{Ciprian Borcea and Ileana Streinu}
\date{}

\begin{document}
\maketitle

\noindent
\begin{abstract}
We study body-and-hinge and panel-and-hinge chains in $R^d$, with two marked points: one on the first body, the other on the last. For a general chain, the squared distance between the marked points gives a Morse-Bott function on a torus configuration space. Maximal configurations, when the distance between the two marked points reaches a global maximum, have particularly simple geometrical characterizations. The three-dimensional case is relevant for applications to robotics and molecular structures.
\end{abstract}



\section*{Introduction}

\medskip \noindent
This work is an extension of our study of singularities of hinge structures \cite{BS}. We refer to that paper
for basic notions and background.

\medskip \noindent
The hinge structures considered here will be body-and-hinge or panel-and-hinge chains in $R^d$ which have a
point marked on the first body and a point marked on the last body. In dimension three, these hinge structures
would model serial  manipulators with revolute joints with a marked base-point and the end-effector as the other
marked point\footnote{The geometrical models have no rotational limitations around the joints and no
self-collision prohibitions.}. Likewise, the panel-and-hinge case may serve as model for ``backbone" protein
chains \cite{BT,CP,BS}. Extremal configurations will be those where the squared distance function between the
two marked points (origin and terminus, or ``head" and ``tail") reaches a local maximum or minimum.

\medskip \noindent
Robotics is obviously concerned with extremal reaches of manipulators and the closely related problem of
identifying the total workspace of a robot. A {\em necessary} condition for extremal configurations was
recognized and proven in several papers \cite{D,KW,SD,S}. In the words of \cite{S}, where the base-point may be
chosen arbitrarily, and the end-point is called ``hand", this necessary condition says: {\em ``the line of sight
from that point to the hand must intersect all turning axes"}\footnote{This incidence of the origin-to-terminus
line with the hinges is understood projectively, that is, includes the possibility of parallelism.}. However,
all critical points with non-zero value for the squared distance function satisfy the condition, and they grow
exponentially with the number of hinges.

\medskip \noindent
In this paper, we refine the study of extremal configurations and obtain, in particular, a very simple {\em
necessary and sufficient} characterization of the global maximum. It may be observed here that the approach used
to identify the global maximum configuration (which is unique for generic body-and-hinge chains) cannot be fully
adapted for the global minimum, although the panel-and-hinge case offers a fair degree of similarity (Theorems
\ref{Maxima} and \ref{Minima}). The distinction, we suggest, stems from the possibility to reinterpret the
global maximum as a global minimum of a related problem. Once recognized, the criterion for the global maximum
can be proven with completely elementary means.

\medskip \noindent
{\bf Global Maximum Theorem:}\ {\em A body-and-hinge chain is in a global maximum configuration if and only if
the segment from the origin to the end-point intersects all hinges in their natural order.}

\medskip \noindent
The panel-and-hinge case has sufficient specificity to warrant separate treatment, most particularly in
dimension three. Local extrema must be global extrema, but are not unique if not flat.

\medskip \noindent
In the final section we discuss some variations.

\section{ A Morse-Bott function for chains with two marked points}

\medskip \noindent
In this section we consider  body-and-hinge chains with   $(n+1)$ bodies and two marked points, one on the first
and the other on the last body. The ambient dimension will be $d$, and we identify the first body with the fixed
reference coordinate system $R^d\equiv B_1$, with the marked point at the origin. The point on the last body
will be the end-point.

\medskip \noindent
The composition of the {\em end-point map} $e: (S^1)^{n}\rightarrow R^d$ with the {\em squared norm} function
$R^d\rightarrow R$ gives the {\em squared distance function} of the end-point to the origin:

\begin{equation}\label{squared-distance}
 F: (S^1)^{n}\rightarrow R, \ \ F(\theta)=<e(\theta),e(\theta)>
\end{equation}

\medskip \noindent  We'll use $T_n=(S^1)^n$ as another notation for the $n$-torus parametrizing the configuration
space of our body-and-hinge chain.

\medskip \noindent
The critical points   of $F$ are described by:

\begin{proposition}
Let $n\geq d$. If {\em zero} is a value of $F$, then all points in $X_0=F^{-1}(0)\subset T_n$ are critical
points of $F$. The critical points with {\em non-zero} critical values are those configurations which have all
hinges (projectively) incident with the line connecting the origin to the end-point.

\medskip \noindent In the generic case,   $F: T_n=(S^1)^n \rightarrow R$ is a Morse-Bott function, which has
only isolated critical points for non-zero critical values, while the fiber over zero $X_0\subset T_n$, when
non-empty, will be smooth, of dimension $n-d$.
\end{proposition}

\medskip \noindent
{\em Proof:}\ The squared norm on $R^d$ has a critical point at the origin, hence all configurations with the
two marked points coinciding (i.e. figuratively, when ``the head bites the tail") will be critical for $F$.

\medskip \noindent
For non-zero critical values, the argument is similar to the one used in \cite{BS} : at a critical
configuration, rotating the part of the chain from hinge $A_i$ on, as a rigid piece, must preserve,
infinitesimally, the (squared) distance ``head-to-tail", that is: must produce a velocity vector for the
end-point orthogonal to the line between the marked points. That requires the (projective) incidence of hinge
and line.

\medskip \noindent
In the generic case, the origin will be a regular value of the end-point map, and, when zero is a value of $F$,
it will give a smooth (not necessarily connected), codimension $d$ fiber $X_0=F^{-1}(0)\subset T_n$. The fact
that the remaining critical points are isolated will follow from the examination of the corresponding Hessian.
The Bott non-degeneracy condition \cite{G} along (the connected components of) $X_0$ will be verified at that
stage as well.

\subsection{Set-up for computing the Hessian matrix}

\medskip \noindent
In a given configuration, hinge number $i$ will be determined by the vector $t_i$ which is the perpendicular
projection of the origin on the corresponding hinge, plus the normal direction to the hyperplane formed by the
hinge and the origin. One may keep track of global orientations, but for our computations, a local choice of
unit normal $\nu_i$ will suffice. Thus $(t_i,\nu_i)=(t_i(\theta),\nu_i(\theta))$ `encodes' the $i^{th}$ hinge.

\medskip \noindent
At any critical point,   we may assume the $\theta\in (S^1)^n$ labelling of the configuration space introduced
by the following rule: the critical point is $\theta=0$, and the position for arbitrary
$\theta=(\theta_1,...,\theta_n)$ is obtained by rotating the last body around the last hinge with angle
$\theta_n$, then rotating the last two bodies (as a rigid piece) around the last but one hinge with angle
$\theta_{n-1}$, and so on until, at last, the whole (rigid) piece thus formed with all the bodies from the
second to the last is rotated around the first hinge with angle $\theta_1$. All rotations, for $i=n,...,1$, are
using the sense dictated by a fixed orientation, say $\{ \nu_i, t_i \}$ in the vector plane $[\nu_i,t_i]$ they
span.

\medskip \noindent
Let $R_i(\omega)$ stand for the {\em linear operator in} $R^d$ which gives the {\em rotation} with angle
$\omega$ around $[\nu_i,t_i]^{\perp}$, with the orientation fixed as above. Then, $\frac{\partial R_i}{\partial
\omega} (\omega)$ is a skew-symmetric operator vanishing on $[\nu_i,t_i]^{\perp}$, and we put $\frac{\partial
R_i}{\partial \omega} (0)=S_i$.

\medskip \noindent
With a second derivation $\frac{\partial^2 R_i}{\partial \omega^2}
(\omega)=-P_iR_i(\omega)$ with $P_i$ denoting the orthogonal
projection on the 2-subspace $[\nu_i,t_i]$. Thus:
$\frac{\partial^2 R_i}{\partial \omega^2}(0)=-P_i$

\medskip \noindent
We let $x\in R^d$ denote the position of the end-point for the critical configuration under investigation. With
the parametrization   and notation just described, and the abbreviation $R_i(\theta_i)=R_i$, the {\em end-point
function} is:

\begin{equation}\label{end-point-function}
 e(\theta)=R_1R_2...R_n x+R_1...R_{n-2}(I-R_n)t_n+...+R_1(I-R_2)t_2+(I-R_1)t_1
\end{equation}

\noindent This gives:

$$ e(0)=x,\ \ \frac{\partial e}{\partial
\theta_i}(0)=S_i(x-t_i) $$

$$ \frac{\partial^2 e}{\partial \theta_i^2}(0)=P_i(t_i-x), \ \ \frac{\partial^2 e}{\partial
\theta_i \partial \theta_j}(0)=S_jS_i(x-t_i), \ j<i $$

\medskip \noindent
Considering that $S_it_i=-||t_i||\nu_i$, the resulting entries for
the Hessian matrix are:

$$ \frac{1}{2}\frac{\partial^2 F}{\partial \theta_i \partial
\theta_j}(0)=<S_i(x-t_i),S_j(x-t_j)>+<S_jS_i(x-t_i),x>=
$$

$$ =<S_i(x-t_i),S_j(x-t_j)>-<S_i(x-t_i),S_jx>=<S_i(t_i-x),S_jt_j>= $$

$$ =<S_i(1-\frac{<x,t_i>}{<t_i,t_i>})t_i,S_jt_j>=(1-\frac{<x,t_i>}{<t_i,t_i>})<S_it_i,S_jt_j>=
$$

$$ =(1-\frac{<x,t_i>}{<t_i,t_i>})||t_i||\cdot ||t_j||
<\nu_i,\nu_j> $$

\noindent for $j\leq i$.

\medskip \noindent
We retain the result of this computation as:

\begin{proposition}\label{Hessian-matrix}
With adequate parametrization, the symmetric $n\times n$ Hesssian matrix for the squared end-point distance
function $F$ at a critical configuration $\theta=0$, with end-point at $e(0)=x$, has entries:

\begin{equation}\label{Hessian}
 \frac{1}{2}\frac{\partial^2 F}{\partial \theta_i \partial
\theta_j}(0)=(1-\frac{<x,t_i>}{<t_i,t_i>})||t_i||\cdot ||t_j|| <\nu_i,\nu_j>, \ \ j\leq i
\end{equation}

\noindent and, after the change of basis $e_i\mapsto
\frac{1}{||t_i||}e_i$, corresponds with the quadratic form given
by:

\begin{equation}\label{Hessian-bis}
 h_{ij}=h_{ji}=(1-\frac{<x,t_i>}{<t_i,t_i>})<\nu_i,\nu_j>, \ \ j\leq i
\end{equation}

\noindent For a non-zero critical value, the coefficients $\alpha_i=\frac{<x,t_i>}{<t_i,t_i>}$ are obtained
geometrically from the intersections of the line through the origin and the end-point $x$ with the hinges, since
these intersections are precisely the points $\frac{1}{\alpha_i}x,\ i=1,...,n$.

\end{proposition}

\medskip \noindent
When $F$ takes the value zero, the Hessian at a critical point in $X_0=F^{-1}(0)$ is equivalent to the Gram
matrix $h_{ij}=<\nu_i,\nu_j>$ of the normals, which is semi-positive definite and has, in the generic case, rank
$d$. Thus, the null-space corresponds to the tangent space of $X_0$ at the critical point under consideration.
This is the Bott non-degeneracy property required along the critical manifold $X_0$.

\subsection{An upper-bound for the number of isolated critical points}

\medskip \noindent
An upper-bound for the number of isolated critical points can be obtained from the following {\em
complexification}: $T_n=(S^1)^n$ is complexified to $(P_1(C))^n$ by considering each circle $S^1$ as the real
locus of the corresponding complex conic:

$$ P_1(C)\approx \{x\in P_2(C)\ :\  x_1^2+x_2^2=x_0^2\} $$

$$ \frac{x_1}{x_0}=\cos \theta_i, \ \ \frac{x_2}{x_0}=\sin \theta_i $$

\noindent With some choice of a reference simplex in each hinge $A_i(\theta)$, say
$a^i_1(\theta),...,a^i_{d-1}(\theta)$, the condition that the `head-to-tail' line meets this hinge becomes

$$ det[a^i_1(\theta) ... a^i_{d-1}(\theta) e(\theta)]=0 $$

\noindent This defines in $(P_1)^n$ a hypersurface of multi-degree $2(d,...,d,1,...1)$ , with positions $d$ up
to the $i^{th}$ coordinate.

\medskip \noindent
The intersection of these $n$ hypersurfaces yields, in the complex count, the number:

$$ 2^n\int_{(P_1)^n} (h_1+...+h_n)(dh_1+h_2+...h_n)...(dh_1+...+dh_{n-1}+h_n)= 2^n\sum_{k=0}^{n-1} A(n,k)d^k $$

\noindent where $h_i$ stands for the class of a point in the $i^{th}$ factor $P_1$, and $A(n,k)$ denote Eulerian
numbers. Thus, we have:

\begin{proposition}
 The number of isolated critical points for
$F:T_n\rightarrow R$ is bounded by

\begin{equation}\label{bound}
 2^n\sum_{k=0}^{n-1} A(n,k)d^k
\end{equation}

\end{proposition}

\subsection{A Meyer-Vietoris sequence (or Morse-Bott theory)}

\medskip \noindent
The determination of critical configurations and their indices can be used in the following setting: we {\em
assume} a generic body-and-hinge chain,   $n\geq 3$ {\em and zero to be a value of $F$}. This ensures the
smoothness of the $(n-3)$-dimensional fiber $X_0=F^{-1}(0)$ over zero, and, for small enough $\epsilon
> 0$, a diffeomorphism $X_0\times B^3_{\epsilon}\approx F^{-1}[0,\epsilon]$, where $B^3_{\epsilon}$ stands for
the 3-dimensional ball of radius $\epsilon$. The $(n-1)$-dimensional fiber $X_{\epsilon}=F^{-1}(\epsilon)$ is
thereby identified with $X_0\times S^2_{\epsilon}$.

\medskip \noindent
We'll use the notation $T_n=(S^1)^n$ for the n-dimensional torus representing the configuration space of our
chain with $(n+1)$ panels and $n$ hinges. For small and nearby values $0 < \gamma < \epsilon < \delta$ we put:

$$ U=F^{-1}[0,\delta),\ V=T_n\setminus F^{-1}[0,\gamma] $$

\noindent assuming all non-zero critical values of $F$ greater than $\delta$. Thus, the two open sets cover the
torus $T_n$, and we have {\em homotopy equivalences}:

$$ U\sim X_0,\ \ U\cap V \sim X_{\epsilon}=X_0\times S^2 $$

\medskip \noindent
The corresponding Meyer-Vietoris exact sequence gives:

$$ ...\rightarrow H_i(X_0\times S^2)\rightarrow H_i(X_0)\oplus H_i(V)\rightarrow H_i(T_n)\rightarrow
H_{i-1}(X_0\times S^2)\rightarrow ... $$

\noindent and in particular, we have the relation of Euler-characteristics:

$$ e(X_0\times S^2)+ e(T_n)=e(X_0)+e(V) \ \ i.e. \ \ e(X_0)=e(V) $$

\medskip \noindent
By Morse theory \cite{M}, the last Euler-characteristic can be expressed in terms of critical points as follows:
put $c_i$ for the number of critical points of index $i$; then:

$$ e(V)=\sum_{i=0}^n (-1)^{n-i}c_i $$

\noindent {\bf Note:}\ The common convention would be to use $-F$ as the Morse function. Since we go with $F$,
the signs are as above.

\medskip \noindent
We obtained:

\begin{proposition}
Suppose $n \geq 3$ and let the origin be a regular value of the end-point map for an otherwise generic chain.
Then, the $(n-3)$-dimensional manifold $X_0$ parametrizing all configurations with the origin coinciding with
the end-point has the Euler number:

\begin{equation}\label{Euler-number}
 e(X_0)=\sum_{i=0}^n (-1)^{n-i}c_i
\end{equation}

\noindent where $c_i$ is the number of critical points of index $i$ for $F$ on $T_n\setminus X_0$.

\end{proposition}

\subsection{The index of the Hessian matrix: panel-and-hinge case}

\medskip \noindent
The restriction to panel-and-hinge chains brings new structural aspects. We note first the presence of natural
transformations of the configurations space $T_n=(S^1)^n$:

\medskip \noindent
{\bf Transformations:}\   Given any configuration, one may consider the hyperplane through one marked point and
some hinge, and reflect the part of the chain from point to hinge in this hyperplane. The resulting
transformation is obviously its own inverse i.e. an {\em involution}. Note that we may always reposition the
structure with its first panel in its fixed location. Since the composition of reflecting in the first and then
last panel gives a global rotation of the chain, these two operations represent the same transformation of $T_n$
and this gives $2n-1$ involutions on the configuration space, all transforming the fibers of $F$ to themselves.
Two such involutions commute when they implicate the same marked point or the respective portions of the chain
do not overlap.

\medskip \noindent
We have seen that, in case $F$ reaches 0, $F^{-1}(0)$ is part of the critical locus, but all critical points for
non-zero critical values are isolated in the generic case. The following definition refers to these isolated
critical configurations corresponding to non-zero critical values of $F$:

\medskip \noindent
\begin{definition}
The $2^n$ flattened configurations   when all panels lie in the same hyperplane (i.e. the codimension one
subspace of the first panel,   which is identified here with $x_d=0$ in $R^d$) will be called   {\em flat
critical configurations (points)}, while critical configurations which do not have all panels in the same
hyperplane will be called {\em non-flat critical configurations (points)}.

\medskip \noindent
It will be observed that, in a non-flat critical configuration, a new panel normal direction requires two
consecutive hinges to meet the line from the origin to the end-point in the same point and we'll call such a
point a {\em fold point}  (on the "head-to-tail" line). Fold points will be considered as labelled by the
corresponding pair of consecutive hinges and ordered via this labelling.
\end{definition}

\medskip \noindent
{\bf Note:}\ Flat critical configurations are fixed points of all involutions described above and the group they
generate, while non-flat critical configurations have  an orbit of cardinality $2^{c}$ under this group, where
$c$ is the number of fold points on the "head-to-tail" line. Assuming a generic case, $c+1$ will give the number
of distinct hyperplanes determined by origin and hinges .

\medskip \noindent
We settle first the case of all flat critical configurations. Then all normal directions are the same and
(regardless of the choice $\pm n_i$) Proposition \ref{Hessian-matrix} gives the signature as that of the
(quadratic form with) matrix:

\begin{equation}\label{flat-Hessian} h_{ij}=h_{ji}=1-\frac{<x,t_i>}{<t_i,t_i>}=\beta_i, \ \ j\leq i
\end{equation}

\begin{lemma}
The signature of the above matrix is determined by the signs of:

$$ \beta_1-\beta_2,\ \beta_2-\beta_3, ..., \beta_{n-1}-\beta_n,
\beta_n $$

\end{lemma}

\medskip \noindent
{\em Proof:}\   The symmetric matrix $H=(h_{ij})$ is the expression of the quadratic form $q(x)=x^tHx$ in the
standard basis $e_i, i=1,...,n$. If we do the (unimodular) change of basis:

$$ e_i\mapsto \tilde{e}_i=e_i-e_{i+1}, \ \ i=1,...,n-1;\ \ e_n=\tilde{e}_n $$

\noindent we obtain a diagonal matrix with the indicated entries.

\medskip \noindent
As mentioned in Proposition 2, if we mark the end-point by $x$, and the intersections of the ``head-to-tail"
line with the $k^{th}$ hinge is $a_kx$, then $\beta_k=1-a_k^{-1}$. For a generic chain in a flattened position,
hinges will intersect the "head-to-tail" line in different points, and we'll have a non-degenerate Hessian.
Since by definition, the index of a symmetric matrix is the number of negative eigen-values, we have:

\begin{proposition}\label{flat-index}
The index of a flat critical configuration   with hinges meeting the line from the origin to the end-point x at
$a_kx,\ k=1,...,n$, is the number of negative values in the list:

\begin{equation}\label{flat-critical-point}
 1-\frac{1}{a_n},\ \frac{1}{a_n}-\frac{1}{a_{n-1}},...,\ \frac{1}{a_2}-\frac{1}{a_1}
\end{equation}

\end{proposition}

\begin{corollary}\label{flat-maximum}
A flat critical configuration is a local maximum   if and only if the hinges meet the oriented segment $(0,x)$
in points $a_kx,\ k=1,...,n$, lined-up in their natural order, that is:

\begin{equation}\label{maximum-ordering} 0< a_1 < a_2 < ... < a_n < 1
\end{equation}

\end{corollary}

\medskip \noindent
For a convenient geometric formulation of the existence of a local minimum in a flat configuration we conceive
of the ``head-to-tail" line as completed to a projective line, and the complement of the affine segment $[0,x]$
gives then the open arc from $0$ to $x$ passing through the ``point at infinity".

\begin{corollary}\label{flat-minimum}
A flat critical configuration is a local minimum   if and only if the hinges meet the oriented arc from 0 to $x$
passing through the point at infinity in their natural order. This means one the following:

\begin{equation}\label{minimum-left}  a_n < ...< a_2< a_1 < 0 \ \ \ \ \mbox{or}
\end{equation}

\begin{equation}\label{minimum-inter}   a_k < ...< a_1 < 0 < 1 < a_n < ... < a_{k+1} \ \ \ \mbox{or}
\end{equation}

\begin{equation}\label{minimum-right}   1< a_n < ... < a_2 < a_1
\end{equation}

\end{corollary}

\medskip \noindent
{\bf Remark:}\ In fact, the projective formulation allows some relaxation in the genericity assumptions and one
of the hinges may be parallel to the ``head-to-tail" line in that flat configuration.

\subsection{Extremal configurations}

\medskip \noindent
In this section we refine our description of extremal configurations of panel-and-hinge chains in $R^d$ with two
marked points.

\medskip \noindent
We'll elaborate on our Corollaries \ref{flat-maximum} and \ref{flat-minimum} and address possible maxima and
minima at non-flat critical configurations. One should remain aware of the involutions described in subsection
1.3. Recall that a fold point on the ``head-to-tail" line is common to two consecutive hinges, say
$a_kx=a_{k+1}x$. When speaking of the ordering of intersections of hinges with the ``head-to-tail" line, either
ordering may be envisaged for the two hinges, but we intend the ordering requested in the statement. The
function $F$ is the squared distance from ``head" to ``tail".

\begin{theorem}\label{Maxima}
A local maximal configuration   for $F$ is characterized by the fact that all hinges intersect the oriented
segment $(0,x)$ in the natural order. Moreover, a local maximum is in fact the global maximum and is unique
modulo the natural transformations generated by the involutions described above. Thus there are $2^{\mu}$
maximal configurations equivalent under natural transformations, where $\mu$ is the number of fold points for
any and all of them.
\end{theorem}

\medskip \noindent
{\em Proof:}\ Let us call {\em flat subsystem}   in a non-flat critical configuration the panel-and-hinge
structure obtained by retaining the consecutive axes which lie in the same hyperplane through the
``head-to-tail" line i.e. the hinges corresponding to a specific normal direction $\nu_i$ (up to the first
occurrence of a fold point, or between a fold point and the next, or up from the last fold point). The ``head"
and ``tail" are inherited from the full configuration. Then the expression (\ref{Hessian-bis}) we obtained for
the Hessian shows that a local maximum requires all flat subsystems to be flat local maxima. Thus, by Corollary
\ref{flat-maximum}, all hinges intersect the oriented segment $(0,x)$ in the natural sequential order.

\medskip \noindent
As a consequence, we may trace a ``red line" on consecutive panels by following the segment $[0,x]$ in our local
maximum configuration. Thus, any other configuration will display the ``red line" as a polygonal arc from
``head" to ``tail" proving that our local maximum is the global maximum. It also follows that any other local
maximum must have exactly the same pattern of planar subsystems and therefore be obtained from our maximal
configuration by some composition of involutions.

\begin{theorem}\label{Minima}
A local minimal configuration   for a non-zero value of $F$ is characterized by the fact that all hinges
intersect the oriented projective arc from 0 to $x$ passing through the point at infinity (i.e. the complement
of $[0,x]$) in the natural order. Moreover, if such a local minimum exists, it is in fact the global minimum and
0 is not a possible value for $F$. When 0 is not a possible value for $F$, all minima are achieved at $2^{\nu}$
minimal configurations equivalent under natural transformations, where $\nu$ is the number of fold points for
any and all of them.
\end{theorem}

\medskip \noindent
{\em Proof:}\ Simple adaptation of the ``red line" argument presented above.

\section{The global maximum as a global minimum}

\medskip \noindent
The criteria obtained above for extremal configurations of marked panel-and-hinge chains in $R^d$ offer obvious
suggestions for the more general {\em body-and-hinge case}.

\medskip \noindent
It is an immediate observation that a marked body-and-hinge chain for which the segment from the origin to the
end-point intersects all hinges in their natural order is in a global maximum configuration, since for any other
configuration, the previous segment (drawn as a ``red segment") becomes a polygonal arc longer than the new
segment from the origin to the end-point. What is less immediate is that any body-and-hinge chain actually
reaches its global maximum in such a configuration. This will be proven by relating the global maximum to a
global minimum.

\begin{theorem}\label{Max-min}
Let a body-and-hinge chain be presented  in a fixed configuration, with the origin as the marked point of the
first panel, $e$ the end-point (on the last body) and hinges given by codimension two affine subspaces $A_i,\
i=1,...,n$. Consider variable points on each hinge $a_i\in A_i$ and the length of the polygonal arc going from
the origin to the end-point through the points $a_i$ in their natural order:

\begin{equation}\label{string}
 f(a_1,...,a_n)=||a_1||+||a_2-a_1||+...+||a_n-a_{n-1}||+||e-a_n||
\end{equation}

\medskip \noindent
The global maximum distance between the marked points of the given body-and-hinge chain equals the global
minimum of the function $f$.

\end{theorem}

\medskip \noindent
{\em Proof:}\ Note that $f$ is continuous and a global minimum always exists.

\medskip \noindent
The proof will follow from the simple case of a single hinge (n=1). In this case, $a_1$ must be the unique point
of $A_1$ which allows a rotation of the segment from $a_1$ to $e$ around $A_1$ to become a continuation of the
segment from the origin to $a_1$. Repeating this observation with respect to $a_{i-1}, a_{i+1}$ and the hinge
$A_i$, shows that the chain can be reconfigured so that the polygonal arc realizing the global minimum becomes a
straight segment from the origin to the end-point intersecting all hinges in their natural order. That is the
global maximum.

\medskip \noindent
We have proven at the same time our:

\medskip \noindent
{\bf Global Maximum Theorem:}\ {\em A body-and-hinge chain is in a global maximum configuration if and only if
the segment from the origin to the end-point intersects all hinges in their natural order. For a generic
body-and-hinge chain, this global maximum is unique.}

\medskip \noindent
{\bf Remarks:}\ (i)\ In general, a body-and-hinge chain may well have local maxima which are not global maxima.
By (\ref{Hessian-bis}), the segment from the origin to the end-point must intersect all hinges.

\medskip
\ (ii)\ $f$ is continuous, but not differentiable when $a_i=a_{i+1}$ for some $i$.

\medskip \noindent
Examples show that global minima for body-and-hinge chains do not have to respect the pattern observed for
panel-and-hinge chains. We may consider, for instance, just two hinges (say in $R^3$) and take the second body
(containing these hinges) as reference. Then the marked points trace circles, each around a hinge. Examples of
global minima with the hinges intersecting in the order $A_2,A_1$ the projective arc from origin to the endpoint
passing thorough infinity, can be easily produced. What still holds true, by (\ref{Hessian-bis}), for any local
minimum is the fact that all hinges intersect the projective arc from the origin to the end-point, passing
through the point at infinity.

\section{Variations on the same theme}

\medskip \noindent
A few `variations' of these techniques should be mentioned before concluding. Again, the issues are important in
robotics and the necessary conditions have been detected in the literature \cite{SR,SD}.

\medskip \noindent
If we abandon the first marked point, we'll rather be concerned with the squared distance from the end-point to
the first hinge. The resulting critical configurations for non-zero values will require the {\em perpendicular}
from the end-point to the first hinge to meet (projectively) all hinges. The global maximum will require the
natural ordering of these intersections on the segment from the foot of the perpendicular to the end-point.

\medskip \noindent
In dimension three, a similar scenario holds if we abandon both marked points and only inquire about the squared
distance between the first hinge and the last hinge. Critical configurations for non-zero values will have the
common perpendicular of these two hinges intersect (projectively) all the intermediate hinges. The global
maximum will require the natural ordering.

\end{document}